\documentclass[11pt,leqno]{amsart}

\usepackage{amssymb,amsmath,amsthm,amsfonts}
\usepackage{latexsym}
\usepackage{color}

\newcommand{\Hmm}[1]{\leavevmode{\marginpar{\tiny%
$\hbox to 0mm{\hspace*{-0.5mm}$\leftarrow$\hss}%
\vcenter{\vrule depth 0.1mm height 0.1mm width \the\marginparwidth}%
\hbox to 0mm{\hss$\rightarrow$\hspace*{-0.5mm}}$\\\relax\raggedright #1}}}
\newcommand{\nc}{\newcommand}
\nc{\les}{\lesssim}
\nc{\nit}{\noindent}
\nc{\nn}{\nonumber}
\nc{\D}{\partial}
\nc{\diff}[2]{\frac{d #1}{d #2}}
\nc{\diffn}[3]{\frac{d^{#3} #1}{d {#2}^{#3}}}
\nc{\pdiff}[2]{\frac{\partial #1}{\partial #2}}
\nc{\pdiffn}[3]{\frac{\partial^{#3} #1}{\partial{#2}^{#3}}}
\nc{\abs}[1] {\lvert #1 \rvert}
\nc{\cAc}{{\cal A}_c}
\nc{\cE}{{\cal E}}
\nc{\cF}{{\mathcal F}}
\nc{\cP}{{\cal P}}
\nc{\cV}{{\cal V}}
\nc{\cQ}{{\cal Q}}
\nc{\cGin}{{\cal G}_{\rm in}}
\nc{\cGout}{{\cal G}_{\rm out}}
\nc{\cO}{{\cal O}}
\nc{\Lav}{{\cal L}_{\rm av}}
\nc{\cL}{{\cal L}}
\nc{\cB}{{\cal B}}
\nc{\cZ}{{\cal Z}}
\nc{\cR}{{\cal R}}
\nc{\cT}{{\cal T}}
\nc{\cY}{{\cal Y}}
\nc{\cX}{{\cal X}}
\nc{\cXT}{{{\cal X}(T)}}
\nc{\cBT}{{{\cal B}(T)}}
\nc{\vD}{{\vec \mathcal{D}}}
\nc{\efield}{\mathcal{E}}
\nc{\vE}{{\vec \efield}}
\nc{\vB}{{\vec \mathcal{B}}}
\nc{\vH}{{\vec \mathcal{H}}}
\nc{\ty}{{\tilde y}}
\nc{\tu}{{\tilde u}}
\nc{\tV}{{\tilde V}}
\nc{\Pc}{{\bf P_c}}
\nc{\bx}{{\bf x}}
\nc{\bX}{{\bf X}}
\nc{\bXYZ}{{\bf XYZ}}
\nc{\bY}{{\bf Y}}
\nc{\bF}{{\bf F}}
\nc{\bS}{{\bf S}}
\nc{\dV}{{\delta V}}
\nc{\dE}{{\delta E}}
\nc{\TT}{{\Theta}}
\nc{\dPsi}{{\delta\Psi}}
\nc{\order}{{\cal O}}
\nc{\Rout}{R_{\rm out}}
\nc{\eplus}{e_+}
\nc{\eminus}{e_-}
\nc{\epm}{e_\pm}
\nc{\eps}{\varepsilon}
\nc{\vnabla}{{\vec\nabla}}
\nc{\G}{\Gamma}
\nc{\w}{\omega}
\nc{\mh}{h}
\nc{\mg}{g}
\nc{\vphi}{\varphi}
\nc{\tlambda}{\tilde\lambda}
\nc{\be}{\begin{equation}}
\nc{\ee}{\end{equation}}
\nc{\ba}{\begin{eqnarray}}
\nc{\ea}{\end{eqnarray}}

\nc{\g}{\gamma}
\nc{\ol}{\overline}

\newtheorem{theorem}{Theorem}[section]
\newtheorem{lemma}[theorem]{Lemma}
\newtheorem{prop}[theorem]{Proposition}

\newtheorem{rmk}[theorem]{Remark}

\def\R{{\rm \rlap{\rm I}\,\bf R}}

\nc{\pT}{\partial_T}
\nc{\pz}{\partial_z}
\nc{\pt}{\partial_t}
\nc{\la}{\langle}
\nc{\ra}{\rangle}
\nc{\infint}{\int_{-\infty}^{\infty}}
\nc{\halfwidth}{6.5cm}
\nc{\figwidth}{10cm}

\nc{\nlayers}{L} \nc{\nsectors}{M}
\nc{\indicator}{\mathbf{1}}
\nc{\Rhole}{R_{\rm hole}}
\nc{\Rring}{R_{\rm ring}}
\nc{\neff}{n_{\rm eff}}
\nc{\Frem}{F_{\rm rem}}
\nc{\Real}{\mathbb R}
\nc{\Z}{\mathbb Z}
\nc{\DD}{\Delta}
\nc{\cD}{\mathcal D}
\nc{\lnorm}{\left\|}
\nc{\rnorm}{\right\|}
\nc{\rnormp}{\right\|_{\ell^{p,\eps}}}
\nc{\rar}{\rightarrow}
\nc{\sgn}{{\rm sign}}
\allowdisplaybreaks
\date{\today}

\begin{document}

\title[High frequency perturbation of cnoidal waves]{High frequency perturbation of cnoidal waves in KdV}

\author{M.~B.~Erdo\smash{\u{g}}an, N.~Tzirakis, and V.~Zharnitsky}
\thanks{The authors were partially supported by NSF grants DMS-0900865 (B.~E.), DMS-0901222 (N.~T.), and DMS-0807897  (V.~Z.)}

\address{Department of Mathematics \\
University of Illinois \\
Urbana, IL 61801, U.S.A.}

\email{berdogan@uiuc.edu \\ tzirakis@math.uiuc.edu\\ vzh@uiuc.edu}

\maketitle

\begin{abstract}
The Korteweg-de Vries (KdV) equation with periodic boundary conditions is considered.
The interaction of a periodic solitary wave (cnoidal wave) with high frequency 
radiation of finite energy ($L^2$-norm) is studied.
It is proved that the interaction of low frequency component (cnoidal wave) and 
high frequency radiation is weak for finite time in the following sense: the 
radiation approximately satisfies Airy equation.

\end{abstract}

\section{Introduction}
The KdV equation
\[
q_t + q_{xxx}+ q_{x}q  = 0
\]
is one of the most basic dispersive partial differential equations (PDE) with solitary wave solutions. 
There are two types
of solitary waves in KdV posed on the real line: exponentially decaying and spatially periodic solitons.
This paper deals exclusively with the periodic case and even more restrictively, we consider
KdV with periodic boundary conditions.
The periodic traveling waves in KdV (already known to Korteweg and de Vries) are called cnoidal waves, as they
may be expressed in terms of the elliptic Jacobi function, see {\em e.g.} \cite{pbs}
\be
\phi_c(z) = \beta_2 +(\beta_3-\beta_2)\, {\rm cn}^2 \left ( \sqrt{\frac{\beta_3-\beta_1}{12}}z;k \right )
\label{eq:cnoidal}
\ee
where
\[
z = x-ct, \beta_1<\beta_2<\beta_3, \beta_1+\beta_2+\beta_3 = 3c, \,\,\,
k^2 = \frac{\beta_3-\beta_2}{\beta_3-\beta_1}.
\]

In this paper we study the dynamics of periodic solutions of the KdV equation
with initial data which is a sum of a  cnoidal wave and of a high frequency "perturbation" with finite energy ($L^2-$ norm).
For small perturbations this falls within the theory of Lyapunov stability results. In the KdV case, this problem was first considered by
Benjamin \cite{benjamin} who conjectured that
cnoidal waves are stable with regard to small perturbations. The conjecture was later proved by McKean \cite{mckean}.
The stability of cnoidal waves on the real line with respect to small perturbations with different periods has also
been considered and is currently an active area of research, see \cite{bottman}
and the references therein.

We, on the other hand, are interested on the behavior of solutions in the case of high frequency
 perturbation. 
 We consider the evolution of solitary
  wave and the high frequency perturbation, and prove that the perturbation evolves almost linearly.
   Within the context of our previous results in \cite{ETZ1, ETZ2} where
we proved that the evolution of high energy solutions of KdV is near-linear, our result in this paper
 can be considered as a superposition principle for a nonlinear dispersive PDE.
It has been long suggested in the physics literature that in 
the regime below collapse, high frequency solutions evolve almost linearly and interact little with
 the low frequencies. This mechanism has been also used, perhaps implicitly, to prove low regularity results,
 see {\em e.g.} \cite{bourgain,ckstt,ckstt1}. In a way, our proof of nearly independent evolution of soliton and
 high frequency radiation, provides some support  of this  heuristic, see e.g. \cite{tao-AMS}, pg. 118.
This nearly independent evolution  is due to  a subtle  averaging effect in the nonlinear dispersive 
dynamics, making these results possible.
Recently KdV was studied with respect to this averaging effect. In \cite{ETZ1, ETZ2},
nearly linear dynamics was established for high frequency initial data and in \cite{titi-arxiv}
 a new elegant proof of well-posedness
 in  $H^s, s\geq 0$ was found using explicitly high frequency averaging effects.

In contrast to the Lyapunov problem, where infinite
time stability is usually established, our result is valid only on finite times,   which is related
to adiabatic invariance phenomena: high frequency wave oscillations are averaged out to produce effective
slow evolution. Since KdV as a model is valid only on a finite time scale,  it is meaningful to consider 
the dynamics of the solutions for finite times.
Given that even the classical adiabatical invariance theorem (conservation of action of the pendulum
with slowly changing frequency) requires careful analysis, our task becomes even harder because of the
infinite dimensionality.

Here, we use these ideas to study the interaction of a solitary wave (low frequencies) with radiation
(high frequencies). Now we state our main theorem:

\begin{theorem}\label{theo:main} Let $\phi(x-ct)$ be a $2\pi$-periodic cnoidal wave solution of  KdV.
Fix $s\in (0,1/2)$. Consider the real valued solution of KdV
on $\mathbf{T}\times \R$ with the initial data $q(x,0)=\phi(x)+g(x)$ satisfying, for some $0<s<1/2$,
$$
\|g\|_2\lesssim 1, \,\,\,\,\,\,\,\|g\|_{H^{-s}}=\eps\ll 1.$$
Then, for each $t>0$, we have
\[
\|q(x, t)-\phi(x-ct)-(e^{tL}g)(x)\|_2 \leq  Ce^{C t} \eps,
\]
where $L=-\partial_x^3- \langle \phi \rangle \partial_x$, and $C$ depends only on $s$ and $\phi$.
\end{theorem}

As usual, $H^{-s}$ is the completion of $L^2$ under the norm $\|u\|_{H^{-s}}=\|\widehat u(k)/(1+|k|^2)^{s/2}\|_{\ell^2}$
and we use the notation   $\langle \phi \rangle = \frac1{2\pi}\int_0^{2\pi} \phi(r)dr  $.
\begin{rmk}
The statement of the theorem above can be extended to an arbitrary $H^4$ solution $\rho$ of KdV in the following sense. Given $T$, there is a constant $C=C(s,T,\rho)$ such that given $\eps>0$ and $g$ as in the theorem we have
\[
\|q(x, t)-\rho(x,t)-(e^{tL}g)(x)\|_2 \leq  C \eps,\,\,\,\,t\in[0,T],
\]
where $L=-\partial_x^3- \langle  \rho \rangle  \partial_x$.

This variation follows from the proof of Theorem~\ref{theo:main} by utilizing Remark~\ref{rmk:K} and Remark~\ref{apriori_rho}.
\end{rmk}

As it is well-known, KdV is a completely integrable system with infinitely many conserved quantities. However, our methods in this paper do not rely on the integrability structure of KdV, and
thus they can be applied to other dispersive models. On the other hand we use the fact that the
smooth solutions of KdV satisfy momentum conservation:
$$\int_{-\pi}^{\pi}u(x,t)dx=\int_{-\pi}^{\pi}u(x,0)dx,$$
and the conservation of energy,
$$\int_{-\pi}^{\pi}u^{2}(x,t)dx=\int_{-\pi}^{\pi}u^{2}(x,0)dx.$$

The KdV equation is locally well-posed in $L^{2}(\Bbb T)$, \cite{Bou2}. Due to energy conservation it is globally well-posed and the solution is in $ C(\Bbb R;L^{2}(\Bbb T)$). Kenig, Ponce, Vega, \cite{kpv}, improved Bourgain's result and showed that the solution of the KdV is locally well-posed in $H^{s}(\Bbb T)$ for any $s>-\frac{1}{2}$. Later, Colliander, Keel, Staffilani, Takaoka, Tao, \cite{ckstt}, showed that the KdV is globally well-posed in
$H^{s}(\Bbb T)$ for any $s \geq -\frac{1}{2}$ thus adding a local well-posedness result for the endpoint $s=-\frac{1}{2}$. Recently T. Kappeler and P. Topalov, \cite{kt} extended the latter result and prove that the KdV is globally well-posed in
$H^{s}(\Bbb T)$ for any $s \geq -1$. Since our statements concerns $L^2$ functions, from the results listed above, we only use the global well-posedness in $L^2$, \cite{Bou2}.

\subsection{Solitary waves description}

Since our goal is to illustrate the phenomena, we restrict ourselves to the special case of
stationary periodic waves with prescribed period $2\pi$. The proof can be readily extended
to arbitrary period.

To make the presentation self-consistent, we directly show that there are such waves, which can also be found
by using properties of Jacobi functions  \eqref{eq:cnoidal}.
Let $q = f(x-ct)$, and substitute this in KdV
\[
-cf^{\prime} = f f^{\prime} + f^{ \prime \prime \prime }.
\]
Integrating once, we obtain
\[
a-cf = \frac{f^2}{2}+ f^{ \prime \prime},
\]
which can be written in the potential form
\be
f^{ \prime \prime} + W_f(f) =0,
\label{eq:2ndorder}
\ee
where
\[
W(f) = \frac{f^3}{6} + c\frac{f^2}{2} - af.
\]
Under the assumption  $c^2+2a > 0$, this cubic polynomial has one local maximum and one local
minimum:
\[
f_- = -\sqrt{c^2+2a}-c, \,\,\, f_+ = \sqrt{c^2+2a}-c.
\]
Taking the 2nd derivative of $W$ at $f_+$,
\[
W_{ff}(f_+) = f_+ + c = \sqrt{c^2+2a},
\]
we obtain the period of small oscillations:
\[
T_0(a,c) = \frac{2\pi}{\sqrt{W_f(f_+)}} =  \frac{2\pi}{\sqrt[4]{c^2+2a}}.
\]
Therefore, moving through the family of periodic solutions nested between the minimum and
the separatrix, we will see the period assuming all intermediate values between  $T_0(a,c)$ and $\infty$.
Thus, if $T_0(a,c) < 1 \Leftrightarrow W^{\prime}(f_+) < 1$ (which can be achieved by taking $a$ or $c$ sufficiently large),
by continuity, somewhere between the critical point and
the separatrix, there will be a $2\pi-$periodic solution. In particular, this can be done by setting $c=0$ and taking $a$ sufficiently large.

We conclude by noting the well known fact that the cnoidal wave $\phi$ is real analytic with exponentially
 decaying Fourier coefficients.

\subsection{Notation}
To avoid the use of multiple constants, we  write $A \lesssim B$ to denote that there is an absolute  constant
$C$ such that $A\leq CB$.  We also write $A\sim B$ to denote both $A\lesssim B$ and $B \lesssim A$.

We define the Fourier sequence of a $2\pi$-periodic $L^2$ function $u$ as
$$u_k=\frac1{2\pi}\int_0^{2\pi} u(x) dx, \,\,\,k\in \mathbb Z.$$
With this normalization we have
$$u(x)=\sum_ke^{ikx}u_k,\,\,\text{ and } (uv)_k=u_k*v_k=\sum_{m+n=k} u_nv_m.$$

\section{Proof of the Main Theorem}\label{sec:evol}
First we discuss that it suffices to prove the theorem for time-independent cnoidal waves.
Consider KdV with periodic boundary conditions (on the circle) $q(x+2\pi)=q(x)$.  Note that if $q(x,t)$ is a solution, then $q(x+c t,t)+c$ is the  solution with initial data $q(x,0)+c$.
In particular, for a cnoidal wave $\phi(x-ct)$, the function $\phi(x)+c$ is also a cnoidal wave. Applying the statement of the theorem with initial data $\phi(x)+c+g(x)$, we obtain at time $t$, that the solution is of the form
$$
q(x+c t,t)+c=\phi(x)+c+(e^{L_ct}g)(x)+O_{L^2}(\eps e^{Ct}),
$$
where $L_c=-\partial_x^3-\big(\frac1{2\pi}\int_0^{2\pi} \phi(r)dr +c  \big)\partial_x$. Noting that
$(e^{L_ct}g)(x)=(e^{Lt}g)(x+ct)$, we obtain
$$
q(x,t)=\phi(x-ct)+(e^{Lt}g)(x)+O_{L^2}(\eps e^{Ct}).
$$
Thus, it suffices to prove the theorem for stationary cnoidal wave.

Let $\phi(x)$ be a time-independent, $2\pi-$periodic cnoidal wave. Consider a solution of KdV of the form   $q(x,t) = \phi(x) + u(x,t)$.  Substituting in KdV, we obtain
\be \label{potential}
u_t + u_{xxx}+ \langle \phi \rangle u_x + (\Phi u)_x +   u u_x =0,\,\,\,\,\,u(x,0)=g(x),
\ee
where $\langle \phi \rangle=\frac1{2\pi} \int_0^{2\pi} \phi(x)dx$, and $\Phi=\phi-\langle \phi \rangle$.

\begin{rmk}

Our assumption on the $H^{-s}$ norm of the initial data $g$, and the momentum conservation implies that   $\langle u(x,t) \rangle = \widehat{g}(0)=O(\eps)$. In the proof of our theorem we will restrict ourselves to the case when $\langle u(x,t) \rangle =0$. This makes the proof more presentable. Removing this assumption introduces more terms in the differentiation by parts formulas which are smaller then the ones we have. In particular, in Theorem~\ref{thm:dbp} below, the formulas for $B(v)$ and $\mathcal R(v)$ would have additional terms which satisfy the a priori estimates given in Proposition~\ref{apriori}.
\end{rmk}

\begin{rmk}
Note that for a mean-zero $L^2$ function $u$, $\|u\|_{H^{-s}}\sim \|u_k/|k|^s\|_{\ell^2}$, we will use this formula without further comments. For a sequence $u_k$, with $u_0=0$, we will use $\|u\|_{H^{-s}}$ notation to denote $\|u_k/|k|^s\|_{\ell^2}$.
\end{rmk}

Using the notation
\[
u(x,t)=\sum_k u_k(t) e^{ikx} \,\,{\rm  and} \,\,\,\, \Phi(x,t)=\sum_k \Phi_k(t) e^{ikx},
\]
we write  \eqref{potential} on the Fourier side,
$$\partial_t u_k=-\frac{ik}{2}   \sum_{k_1+k_2=k}u_{k_1}u_{k_2}- ik    \sum_{k_1+k_2=k}\Phi_{k_1}u_{k_2}+ i(k^3-ak)u_k,\,\,\,\,\,\,u_k(0)=\widehat g(k),$$
where $a=\langle \phi \rangle$. Because of the mean zero assumption on $\Phi$ and $u$, and conservation of momentum, there are no zero harmonics in this equation. Without the mean
zero assumption this equation would have an additional term of the form $iu_0k\Phi_k$ which is of order
$\eps$ and has fast decay in $k$.

Using the transformations
\begin{align*}
u_k(t)&=v_k(t)e^{i(k^3-ak)t},\\
\Phi_k(t)&=S_k(t)e^{i(k^3-ak)t},
\end{align*}
and the identity
$$(k_1+k_2)^3-k_1^3-k_2^3-a(k_1+k_2)+ak_1+ak_2=3(k_1+k_2)k_1k_2,$$
the equation can be written in the form
\be\label{v_eq}
\partial_t v_k=-\frac{ik}{2}   \sum_{k_1+k_2=k}e^{-i3kk_1k_2t}(v_{k_1}+2S_{k_1})v_{k_2}.
\ee
The following theorem  will be proved in Section~\ref{sec_dbp} by distinguishing the resonant and nonresonant sets and using differentiation by parts.

\begin{theorem}\label{thm:dbp}
The system \eqref{v_eq} can be written in the following form:
\be\label{v_eq_dbp}
\partial_t[v+K(v)+B(v)]_k=L_0(v)_k+\mathcal{R}(v)_k,
\ee
where we define $K(v)_0=B(v)_0=L_0(v)_0+\mathcal{R}(v)_0=0$, and for $k\neq 0$, we define
\begin{align*}
K(v)_k& = - \sum_{k_1+k_2=k}\frac{e^{-3ikk_1k_2t}  S_{k_1} v_{k_2}}{3k_1k_2},
\\
B(v)_k&= - \sum_{k_1+k_2=k}\frac{e^{-3ikk_1k_2t}  v_{k_1} v_{k_2}}{6k_1k_2}
\\
& - \frac19\sum_{k_1+k_2+k_3=k}^{*}\frac{e^{-3it(k_1+k_2)(k_2+k_3)(k_3+k_1)}}{k_1(k_1+k_2)(k_2+k_3)(k_3+k_1)}(v_{k_1}+S_{k_1})v_{k_2}v_{k_3},
\\
L_0(v)_k&=\frac{2i}{3}\sum_{k_1+k_2+k_3=k}\frac{e^{-3it(k_1+k_2)(k_2+k_3)(k_3+k_1)}}{k_1} S_{k_1}S_{k_2}v_{k_3},
\end{align*}
\begin{align*}
&\mathcal{R}(v)_k=\frac{i}{3}\frac{v_{k}|v_{k}|^2}{k}-\frac{i}{3}\frac{S_{-k}v_{k}v_{k}}{k} - \frac{2i}{3}v_k\sum_{ |j|\neq k} \frac{S_{j}v_{-j}}{j}-\frac{2i}{3}\frac{S_k}{k}\sum_j S_{-j}v_j\\
&+\sum_{k_1+k_2=k} \frac{e^{-3ikk_1k_2t}}{3k_1k_2}(\partial_{t} S_{k_1}) v_{k_2}
\\
&+\frac{2i}{3}\sum_{\stackrel{k_1+k_2+k_3=k}{k_2+k_3\neq 0}} \frac{e^{-3it(k_1+k_2)(k_2+k_3)(k_3+k_1)}}{k_1} v_{k_1}S_{k_2}v_{k_3}
\\
&-\frac{1}{9}
\sum_{k_1+k_2+k_3=k}^{*}\frac{e^{-3it(k_1+k_2)(k_2+k_3)(k_3+k_1)}v_{k_2}v_{k_3}\partial_{t}S_{k_1}}{k_1(k_1+k_2)(k_2+k_3)(k_3+k_1)}
\\
&+\frac{i}{9} \sum_{k_1+k_2+k_3+k_4=k}^{\star}\frac{e^{it\tilde\psi(k_1,k_2,k_3,k_4)}(k_3+k_4)S_{k_1}v_{k_2}(v_{k_3}+2S_{k_3})v_{k_4}}{k_1(k_1+k_2)(k_1+k_3+k_4)(k_2+k_3+k_4)}
\\
&+\frac{i}{18} \sum_{k_1+k_2+k_3+k_4=k}^{\star}\frac{e^{it\tilde\psi(k_1,k_2,k_3,k_4)}(k_1+2k_3+2k_4)v_{k_1}v_{k_2}(v_{k_3}+2S_{k_3})v_{k_4}}{k_1(k_1+k_2)(k_1+k_3+k_4)(k_2+k_3+k_4)}.
\end{align*}
Here $\sum^*$ means the sum does not contain the terms which makes the denominator zero.
\end{theorem}

\begin{prop}\label{apriori}
Assume that $\|v\|_2\lesssim 1 $ and $0<s<1/2$, then
\begin{align}\label{p1}
&\|K(v)\|_{H^{-s}}\leq \|K(v)\|_2 \lesssim \|v\|_{H^{-s}},
\\
\label{p2}
&\|B(v)\|_{H^{-s}}\leq \|B(v)\|_2 \lesssim \|v\|_{H^{-s}}^2,
\\\label{p3}
&\|L_0(v)\|_{H^{-s}} \lesssim \|v\|_{H^{-s}},\,\,\,\,\,\,\,\|L_0(v)\|_{2} \lesssim \|v\|_{2},
\\\label{p4}
&\|\mathcal{R}(v)\|_{H^{-s}}\leq \|\mathcal{R}(v)\|_{2} \lesssim \|v\|_{H^{-s}}.
\end{align}
\end{prop}
We will prove this proposition in Section~\ref{sec_apriori}. Now, we continue with the proof of the main theorem. First we will prove the near-linear behavior using a modified linear operator (Theorem~\ref{thm:L1}),  then we will prove that the modified linear evolution is close to the Airy evolution (Theorem~\ref{thm:airy}). These two theorems imply Theorem~\ref{theo:main}.
\begin{theorem}\label{thm:L1} Let $0<s<1/2$.
Let $u$ be a mean zero solution of \eqref{potential} with $u(\cdot,0)=g$, where
$$\|g\|_{L^2}\lesssim 1, \text{ and } \,\,
\|g\|_{H^{-s}}\lesssim \eps\ll 1.$$
Then, for each $t>0$, we have
\[
\|u(\cdot, t)-e^{tL_1}g\|_2 \leq  Ce^{C t} \eps.
\]
Here $L_1=L+P$, where $L=-\partial_x^3-\langle \phi(\cdot)\rangle \partial_x$ and $P$ is defined in the Fourier side as
$$
(Pu)_k = \big [ e^{Lt} L_0(e^{-Lt} u) \big ]_k =   \frac{2i}{3} \sum_{k_1+k_2+k_3=k}\frac{\Phi_{k_1}}{k_1}\Phi_{k_2}u_{k_3},\,\,\,k\neq 0,
$$
and $(Pu)_0=0.$
\end{theorem}
\begin{theorem}\label{thm:airy} Let $g$, $L_1$, $L$ be as in the previous theorem.
Then, for each $t>0$, we have
\[
\|e^{tL_1}g-e^{tL}g\|_2 \leq  Ce^{C t} \eps.
\]
\end{theorem}
\begin{proof}[Proof of Theorem~\ref{thm:L1}]
First we will prove that the norm assumptions on the initial data remain  intact up to times of order $\log(1/\eps)$.
By  $L^2$ conservation in KdV, we have
$$\|u(\cdot,t)\|_{L^2}\leq\|q(\cdot,t)\|_{L^2}+\|\phi(\cdot)\|_{L^2}=\|q(\cdot,0)\|_{L^2}+\|\phi(\cdot)\|_{L^2}\lesssim 1.$$
Now we prove that for some $C$,
\be\label{hsbound}
\|u(\cdot,t)\|_{H^{-s}}\leq Ce^{Ct} \eps.
\ee
To prove this integrate \eqref{v_eq_dbp} from $0$ to $T$ to obtain
\be\nn
(I+K)v(T)=(I+K)v(0)+B(v)(0)-B(v)(T)+  \int_0^T (L_0(v)+\mathcal{R}(v)) dt.
\ee
Since $u_k=v_k e^{i\psi(k)t}$ and $\Phi_k=S_ke^{i\psi(k)t}$ with $\psi(k)=k^3-ak$, we have
\begin{align}\label{u_eqn_s}
(I+\tilde K)u(T)&=e^{i\psi(k)T}(I+K)u(0)+e^{i\psi(k)T}B(u)(0)\\
&-e^{i\psi(k)T}B(e^{-i\psi(k)T} u )(T) \nonumber\\
&+ e^{i\psi(k)T} \int_0^T (L_0(e^{-i\psi(k)t} u )+\mathcal{R}(e^{-i\psi(k)t} u ) dt, \nonumber
\end{align}
where $\tilde K$ is the time-independent operator:
$$
\tilde K (u)_k = -\sum_{k_1+k_2=k} \frac{e^{-3ik_1k_2kT-i\psi(k_1)T-i\psi(k_2)T+i\psi(k)T}\Phi_{k_1}u_{k_2}}{3k_1k_2}
= -\sum_{k_1+k_2=k} \frac{\Phi_{k_1}u_{k_2}}{3k_1k_2}.
$$

\begin{lemma} \label{fred} For $0<s\leq 1$,
$$\|(I+\tilde K)u\|_{H^{-s}}\gtrsim \|u\|_{H^{-s}}.$$
\end{lemma}
\begin{proof}
This follows from Fredholm alternative. First note that for $u\in L^2$ with mean zero
$$\|\tilde K u\|_2= \Big\|\sum_{k_1+k_2=k} \frac{\Phi_{k_1}u_{k_2}}{3k_1k_2}\Big\|_2 \leq \Big\|\frac{\Phi_{k} }{3k }\Big\|_{\ell^1}\Big\|\frac{ u_{k}}{ k}\Big\|_2 \lesssim \|u\|_{H^{-s}}.
$$
By the density of $L^2$ in $H^{-s}$, this inequality holds for each $u\in H^{-s}$.
Therefore, by Rellich's theorem, $\tilde K$ is a compact operator on $H^{-s}$.
It suffices to show that the kernel of $I+\tilde K$ is trivial. Note that if $(I+\tilde K)u=0$ for some $u\in H^{-s}$, then by the discussion above,  $\tilde K u\in L^2$, and hence $u\in L^2$.
Using the definition of $\tilde K$, we have
\[
(I+\tilde K)u =0 \Leftrightarrow u(x) - \frac{1}{3} \Phi_{-1}(x) \, u_{-1}(x) = 0,
\]
where $f_{-1}(x)$ denotes the mean-zero antiderivative of a mean-zero function, 
\be\label{avzero}
f_{-1}(x)=\int_0^x f(r) dr+\frac{1}{2\pi}\int_0^{2\pi} r f(r) dr.
\ee

 Let $U(x)=u_{-1}(x)$, then we obtain equivalent 1st order linear ODE
\[
U^{\prime} -  \frac{1}{3} \Phi_{-1}(x) U(x) = 0,
\]
which has the general solution $U(x) = U(0) \exp(\int_0^x \Phi_{-1}(x) dx)$, which
can be mean-zero only if $U(x) \equiv 0$. This implies $u(x) \equiv 0$.

\end{proof}
\begin{rmk}\label{rmk:K}
Let
$$\tilde K_t(u):= -\sum_{k_1+k_2=k}\frac{\rho_{k_1}(t)u_{k_2}}{3k_1k_2},$$
where $\rho_k(t)=\widehat{\rho(\cdot,t)}(k)$, and $\rho$ is an $L^2$ solution of KdV.
Then, the statement of the lemma is valid for $\tilde K_t$\\
i) for  $t\in [0,T]$ with a constant $C_T$ depending on $T$,\\
ii) on $\mathbb R$ with a constant independent of $t$ if
$\|\rho(\cdot,t)\|_2=\|\rho(\cdot,0)\|_2 <1/2$.\\
To prove these statements first note that
\be \label{tildeK}
\|\tilde K_t\|_{H^{-s}\to H^{-s}} \leq 2 \|\rho(\cdot, t)\|_2.
\ee
Then, using the resolvent identity,
$$(I+\tilde K_t)^{-1}-(I+\tilde K_s)^{-1}=(I+\tilde K_s)^{-1}(\tilde K_s-\tilde K_t)(I+\tilde K_t)^{-1},$$
the linearity of $\tilde K_t$ in $\rho$,  and \eqref{tildeK}, we see that  the operator $(I+\tilde K_t)^{-1}$ is continuous in time in the operator norm $H^{-s}\to H^{-s}$. Thus, $(i)$ follows form the proof of the lemma and compactness.   To see (ii), note that by \eqref{tildeK}  the operator norm of $\tilde K$ in $H^{-s}$ is $<1$, and invert the operator using Neumann series.
\end{rmk}

Using Lemma~\ref{fred} and Proposition~\ref{apriori} in \eqref{u_eqn_s}, we obtain
\begin{align*}
\|u(T)\|_{H^{-s}}&\lesssim \|u(0)\|_{H^{-s}}+\|u(0)\|_{H^{-s}}^2+\|u(T)\|_{H^{-s}}^2
+\int_0^T \|u(t)\|_{H^{-s}} dt\\
&\leq C \eps +C \|u(T)\|_{H^{-s}}^2+C \int_0^T \|u(t)\|_{H^{-s}} dt.
\end{align*}
Therefore on $[0,T_0]$, where $T_0=\inf \{T:\|u(t)\|_{H^{-s}}\geq \frac{1}{2C}\}$, we have
$$
\|u(T)\|_{H^{-s}}\leq  2C \eps  + 2C \int_0^T \|u(t)\|_{H^{-s}} dt.
$$
This implies by Gronwall that $\|u(T)\|_{H^{-s}}\leq 2C\eps e^{2CT}$ as claimed.

Now, note that $\tilde K(u)=e^{-i\psi(k)T}K(e^{i\psi(k)t}u)$. Therefore, using  Proposition~\ref{apriori} (for $K$, $B$, and $\mathcal R$) and \eqref{hsbound} in \eqref{u_eqn_s}, we obtain
\be\label{newu}
u(T)=e^{i\psi(k)T}u(0) + e^{i\psi(k)T} \int_0^T L_0(e^{-i\psi(k)t} u ) dt +O_{L^2}(\eps e^{CT}).
\ee
Using the definition of $L_0$, we have
\begin{align*}
&L_0(e^{-i\psi(k)t} u ) =
\\
&=\frac{2i}{3} \sum_{k_1+k_2+k_3=k}\frac{e^{it[-3 (k_1+k_2)(k_2+k_3)(k_3+k_1)-\psi(k_1)-\psi(k_2)-\psi(k_3)]}}{k_1} \Phi_{k_1}\Phi_{k_2}u_{k_3}\\
&=\frac{2i}{3}e^{-i\psi(k)t}\sum_{k_1+k_2+k_3=k}\frac{\Phi_{k_1}\Phi_{k_2}u_{k_3} }{k_1}=e^{-Lt}P u.
\end{align*}
In the last line we used the fact that for $k=k_1+k_2+k_3$,
$$
-\psi(k) = -3 (k_1+k_2)(k_2+k_3)(k_3+k_1)-\psi(k_1)-\psi(k_2)-\psi(k_3),
$$
and the definition of $L$ and $P$.

Using this in \eqref{newu}, we have
\begin{align*}
u(T)=&e^{LT}u(0) + e^{LT} \int_0^T e^{-Lt} P(u(t)) dt +O_{L^2}(\eps e^{CT})
\\
=&e^{LT}u(0) + \int_0^T e^{L(T-t)} P(e^{L_1t}u(0)) dt
\\
&+  \int_0^T e^{L(T-t)} P[u(t)-e^{L_1t}u(0)] dt+O_{L^2}(\eps e^{CT})
\\
=&e^{L_1T}u(0)+  \int_0^T e^{L(T-t)} P[u(t)-e^{L_1t}u(0)] dt+O_{L^2}(\eps e^{CT}).
\end{align*}
Let $h(t):=\|u(t)-e^{L_1t}u(0)\|_2$. Using the equality above and the bound for $L_0$ in Proposition~\ref{apriori}  for the operator $P(u(t))=e^{Lt}L_0(e^{-Lt}u(t))$, we obtain
$$
h(T)\lesssim \eps e^{CT} +\int_0^T h(t)dt.
$$
The theorem follows from this by Gronwall.
\end{proof}

\begin{proof}[Proof of Theorem~\ref{thm:airy}]
First we prove that our assumptions on the initial data remain intact for times of order $\log(1/\eps)$:
\begin{lemma}\label{interpol}
For $s\in[-1,1]$, the operator $L_1$ defined above satisfies
$$
\|e^{tL_1}\|_{H^{s}\to H^{s}} \lesssim e^{Ct}.
$$
\end{lemma}
\begin{proof}
First note that we can rewrite $ P(u)$ (for mean zero $u$ and $\Phi$) in the following form which is valid for each $k\in \mathbb Z$,
\begin{align*}
P(u)_k & = -\frac{2}{3} \sum_{k_1+k_2+k_3=k}\frac{\Phi_{k_1}}{ik_1}\Phi_{k_2}u_{k_3}+\frac{2}{3} \sum_{k_1+k_2+k_3=0} \frac{\Phi_{k_1}}{ik_1}\Phi_{k_2}u_{k_3}.
 \end{align*}
The constant term  makes the right-hand side vanish for $k=0$, which makes $P(u)$ mean-zero in the space side. Using the formula \eqref{avzero} for the function $\Phi$, we write $L_1=L+P$ in the space side as
$$
L_1u =-\partial_x^3 u - a\partial_x u + G u -\frac1{2\pi}\langle u, G\rangle,
$$
where
\begin{align*}
G(x)=-\frac23\Phi(x)\Big(\int_0^x\Phi(r)dr+\frac1{2\pi}\int_0^{2\pi}r\Phi(r)dr\Big)
\end{align*}
Also note that
$$
L_1^*u =\partial_x^3 u + a\partial_x u + G u -\frac1{2\pi}G \langle u, 1\rangle.
$$
Note that by duality and interpolation it suffices to prove the assertion of the lemma for $s=1$ for $L_1$ and $L_1^*$.
We will give the proof for $L_1$, for $L_1^*$ the proof is essentially the same since they have very similar forms.
Consider the equation
$$
u_t=L_1u=-\partial_x^3 u - a\partial_x u + G  u -\frac1{2\pi}\langle u, G\rangle.
$$
First, we calculate
\begin{align*}
&\frac12\frac{d}{dt}\|u\|_2^2  = \int_0^{2\pi} u_{t}u  dx \\
&= - \int_0^{2\pi} u_{xxx} u  dx - a\int_0^{2\pi} u_{x} u dx + \int_0^{2\pi}G u^2  dx-\frac1{2\pi}\langle u, G\rangle \int_0^{2\pi} u dx\\
&= \int_0^{2\pi}G u^2  dx-\frac1{2\pi}\langle u, G\rangle \int_0^{2\pi} u dx.
\end{align*}
This implies that
$$
\Big|\frac{d}{dt}\|u\|_2^2\Big|\leq 2\|G\|_\infty \|u\|_2^2+\|G\|_2  \|u\|_2^2 \lesssim \|u\|_2^2,
$$
Similarly,
\begin{align*}
&\frac12\frac{d}{dt}\|u_x\|_2^2 =  \int_0^{2\pi} u_{xt}u_x dx \\
&= - \int_0^{2\pi} u_{xxxx} u_x  dx - a\int_0^{2\pi} u_{xx} u_x dx + \int_0^{2\pi}(Gu)_x u_x dx \\
&=\int_0^{2\pi}G u_x^2 dx-\frac12\int_0^{2\pi}u^2 G_{xx} dx \\
&= O\Big(\|u_x\|_2^2 \|G \|_\infty+\|u\|_2^2 \| G_{xx}\|_\infty  \Big) = O\big(\|u\|_{H_1}^2\big).
\end{align*}
Combining the two inequalities, we obtain
$$
\Big|\frac{d}{dt}\|u\|_{H^1}^2\Big|\lesssim \|u\|_{H^1}^2,
$$
which finishes the proof.
\end{proof}

We return to the proof of Theorem~\ref{thm:airy}.
Consider the equation
$$u_t=L_1 u,\,\,\,\,\,\,u(0,x)=g(x).$$
Repeating the discussion in the beginning of Section~\ref{sec:evol}, and introducing the variables   $v_k$ and $S_k$ as above, we have
\be\label{L1v}
\partial_t v_k=\frac{2i}{3}\sum_{k_1+k_2+k_3=k}e^{-3it(k_1+k_2)(k_2+k_3)(k_1+k_3)}\frac{S_{k_1}}{k_1}S_{k_2}v_{k_3}.
\ee
We will prove the following proposition, using differentiation by parts, in Section~\ref{sec_dbp}.
\begin{prop}\label{dbpL1}
The system \eqref{L1v} is equivalent to the following
\be\label{L1v2}
\partial_t\big(v_k+D(v)_k\big)=E(v)_k,
\ee
where
$$
D(v)_k=-\frac{i}{3}\sum_{k_1+k_2+k_3=k}^{*}\frac{e^{-3it(k_1+k_2)(k_2+k_3)(k_3+k_1)}}{k_1(k_1+k_2)(k_2+k_3)(k_3+k_1)} S_{k_1} S_{k_2}v_{k_3},
$$

\begin{align*}
E(v)_k&=\frac{2i}{3}S_k\sum_{j\neq k} \frac{S_jv_{-j}}{j}+\frac{2i}{3}\frac{S_k}{k}\sum_{j}  S_{-j}v_{j}\\
&-\frac{i}{3}\sum_{k_1+k_2+k_3=k}^{*}\frac{e^{-3it(k_1+k_2)(k_2+k_3)(k_3+k_1)}\partial_t(S_{k_1} S_{k_2})v_{k_3}}{k_1(k_1+k_2)(k_2+k_3)(k_3+k_1)}
\\&+\frac{2}{9}\sum_{k_1+k_2+k_3+k_4+k_5=k}^{*}\frac{e^{ it\tilde\psi(k_1,k_2,k_3,k_4,k_5)}S_{k_1} S_{k_2} S_{k_3}S_{k_4}v_{k_5}}{k_1k_3(k_1+k_2)(k-k_1)(k-k_2)}.
\end{align*}
Here $\tilde \psi$ is a real valued phase function which is irrelevant for the proof of the Theorem.
\end{prop}
\begin{prop}\label{apriori2}
The following a priori estimates hold
$$
\|D(v)\|_2\lesssim \|v\|_{H^{-1}},\,\,\,\,\,\, \|E(v)\|_2\lesssim \|v\|_{H^{-1}}.
$$
\end{prop}
\begin{proof}
Using $|k_1||k_1+k_3|\geq |k_3|$ (for nonzero integer values of $k_1$ and $k_1+k_3$), and eliminating the product $|(k_1+k_2)(k_2+k_3)|$ from the denominator, we have
$$
|D(v)_k| \lesssim \sum_{k_1+k_2+k_3=k} |S_{k_1}S_{k_2}|\frac{|v_{k_3}|}{|k_3|}.
$$
This implies by Young's inequality
$$
\|D(v)\|_2\lesssim \|S\|_{\ell^1}\|S\|_{\ell^1} \|v_k/k\|_2\lesssim \|v\|_{H^{-1}}.
$$
The proof for the contribution of the second line in the definition of $E(v)$ is exactly the same using the fast decay of $|\partial_tS_k|\lesssim |k^3| |S_k|$.
The $L^2$ norm of the   first line is
\begin{align*}
&\leq \|S\|_2 \sum_j \frac{|S_j||v_j|}{|j|}+\|S_k/k\|_2\sum_j |jS_j|\frac{|v_j|}{|j|}  \\
& \lesssim [\|S\|_2^2 + \|S_k/k\|_2\|kS_k\|_2]\|v_k/k\|_2\lesssim \|v\|_{H^{-1}}.
\end{align*}
The second inequality follows from Cauchy Schwarz.
To estimate the third line in the definition of $E(v)$, note that
$$
|k-k_2||k_1|k_3|=\frac{|k_1+k_3+k_4+k_5||k_1| |k_3| |k_4|}{|k_4|}\gtrsim\frac{|k_5|}{|k_4|}
$$
for nonzero integer values of the factors.
Using this in the sum and eliminating the rest of the factors, we estimate the third line as
$$
\lesssim \sum_{k_1+k_2+k_3+k_4+k_5=k}  |S_{k_1} S_{k_2} S_{k_3} k_4S_{k_4}|\frac{|v_{k_5}|}{|k_5|}.
$$
As above, this implies that the $L^2$ norm is
$$
\lesssim \|S\|^3_{\ell^1} \|kS_k\|_{\ell^1} \|v_k/k\|_2\lesssim \|v\|_{H^{-1}}.\qedhere
$$
\end{proof}
To complete the proof of Theorem~\ref{thm:airy}, integrate \eqref{L1v2} from $0$ to $T$:
$$
v_k(T)=v_k(0)-D(v)_k(T)+D(v)_k(0)+\int_0^T E(v)_k(t) dt.
$$
Using the transformation $u_k=v_k e^{i\psi(k)t}$, we have
\begin{align*}
u_k(T)&=e^{i\psi(k)T}u_k(0)-e^{i\psi(k)T}D(e^{-i\psi(k)T}u)_k(T)+e^{i\psi(k)T}D( u)_k(0)\\
&+\int_0^T e^{i\psi(k)T}E(e^{-i\psi(k)t}u)_k(t) dt.
\end{align*}
Noting that $(e^{LT}u)_k=e^{i\psi(k)T}u_k(0)$, and using the a priori estimates in Proposition~\ref{apriori2}, we have
\begin{align*}
\|u(T)-e^{LT}u(0)\|_2&\lesssim \|u(T)\|_{H^{-1}}+\|u(0)\|_{H^{-1}}+\int_0^T \|u(t)\|_{H^{-1}} dt\\
&\lesssim \eps e^{CT}.
\end{align*}
In the last line we used Lemma~\ref{interpol} and the hypothesis on $\|u(0)\|_{H^{-s}} \geq \|u(0)\|_{H^{-1}}$.
\end{proof}

\section{Differentiation by Parts} \label{sec_dbp}
In this section we prove Theorem~\ref{thm:dbp} and Proposition~\ref{dbpL1}.
\begin{proof}[Proof of Theorem~\ref{thm:dbp}]
We need to obtain \eqref{v_eq_dbp} from \eqref{v_eq} by using differentiation by parts.
It will be  useful to name the terms appearing in the formula \eqref{v_eq_dbp}.
We will denote the terms in the first and second line of the definition of $B(v)$ by $B_1(v)$ and $B_2(v)$, respectively. We also denote the term in the $j$th line of the definition of $\mathcal R(v)$ in the statement of the Theorem by $\mathcal R_j (v)$, $j=1,2,...,6$, and we further denote the four summands in $\mathcal R_1(v)$ by $\mathcal R_{1,m} (v)$, $m=1,2,3,4$.

Since $e^{-3ikk_1k_2t}=\partial_{t}( \frac{i}{3kk_1k_2}e^{-i3kk_1k_2t})$, using differentiation by parts we can rewrite \eqref{v_eq} as
\begin{align*}
\partial_{t}v_{k}&=\partial_{t}\Big( \frac12  k \sum_{k_1+k_2=k}\frac{e^{-3ikk_1k_2t}(v_{k_1}+2S_{k_1})v_{k_2}}{3 kk_1k_2}\Big)\\
& - \frac12 k\sum_{k_1+k_2=k}\frac{e^{-3ikk_1k_2t}}{3kk_1k_2}\partial_{t}[(v_{k_1}+2S_{k_1}) v_{k_2}].
\end{align*}
Recalling the definition of $K(v)$ and $B_1(v)$ from equation \eqref{v_eq_dbp}, we can rewrite this equation in the form:
$$
\partial_{t}[v_{k}+K(v)_k+B_1(v)_k]=\sum_{k_1+k_2=k}\frac{e^{-3ikk_1k_2t}}{6k_1k_2}\partial_{t}[(v_{k_1}+2S_{k_1}) v_{k_2}].
$$
Note that since $v_0=0$ and $S_0=0$, in the sums above $k_1$ and $k_2$ are not zero. We now handle the term  when the derivative hits $v_{k_1}v_{k_2}$. By symmetry and \eqref{v_eq}, we have
\begin{align*}
&\sum_{k_1+k_2=k}\frac{e^{-3ikk_1k_2t}}{k_1k_2}\partial_{t}(v_{k_1}v_{k_2})= \\ &= -i\sum_{k=k_1+k_2}\frac{e^{-3ikk_1k_2t}}{k_1}v_{k_1}\Big( \sum _{\mu+\lambda=k_2}e^{-3itk_2\mu\lambda}(v_{\mu}+2S_\mu)v_{\lambda}\Big)\\
&= -i\sum_{k=k_1+\mu+\lambda}\frac{v_{k_1}(v_{\mu}+2S_\mu) v_{\lambda}}{k_1}e^{-3it[kk_1(\mu+\lambda)+\mu\lambda(\mu+\lambda)]}.
\end{align*}
We note that $\mu+\lambda$ can not be zero since $\mu+\lambda=k_2$. Using the identity
$$kk_1+\mu\lambda=(k_1+\mu+\lambda)k_1+\mu\lambda=(k_1+\mu)(k_1+\lambda)$$
and  by renaming the variables $k_2=\mu, k_3=\lambda$, we have that
\begin{align*}
&\sum_{k_1+k_2=k}\frac{e^{-3ikk_1k_2t}}{k_1k_2}\partial_{t}(v_{k_1} v_{k_2})=\\&=-i\sum_{\stackrel{k_1+k_2+k_3=k}{k_2+k_3\neq 0}}\frac{e^{-3it(k_1+k_2)(k_2+k_3)(k_3+k_1)}}{k_1}v_{k_1}(v_{k_2}+2S_{k_2})v_{k_3}.
\end{align*}
Calculating the  term when the derivative hits $S_{k_1}v_{k_2}$ similarly, we have
\be\label{y1y5}
\partial_t[v_{k}+K(v)_k+B_1(v)_k]= \sum_{j=1}^5Y_j(v)_k,
\ee
where
$$Y_1(v)_k=-\frac{i}{3}\sum_{\stackrel{k_1+k_2+k_3=k}{k_2+k_3\neq 0}}\frac{e^{-3it(k_1+k_2)(k_2+k_3)(k_3+k_1)}}{k_1} v_{k_1}v_{k_2}v_{k_3},$$
$$Y_2(v)_k=-\frac{i}{3}\sum_{\stackrel{k_1+k_2+k_3=k}{k_2+k_3\neq 0}}\frac{e^{-3it(k_1+k_2)(k_2+k_3)(k_3+k_1)}}{k_1} S_{k_1}v_{k_2}v_{k_3},$$
$$Y_3(v)_k=-\frac{2i}{3}\sum_{\stackrel{k_1+k_2+k_3=k}{k_2+k_3\neq 0}}\frac{e^{-3it(k_1+k_2)(k_2+k_3)(k_3+k_1)}}{k_1} v_{k_1}S_{k_2}v_{k_3},$$
$$Y_4(v)_k=-\frac{2i}{3}\sum_{\stackrel{k_1+k_2+k_3=k}{k_2+k_3\neq 0}}\frac{e^{-3it(k_1+k_2)(k_2+k_3)(k_3+k_1)}}{k_1} S_{k_1}S_{k_2}v_{k_3},$$
$$Y_5(v)_k=-\sum_{k_1+k_2=k}\frac{e^{-3ikk_1k_2t}}{3k_1k_2}(\partial_{t} S_{k_1}) v_{k_2}.$$
Note that
\be\label{y3y4y5}
Y_3(v)=\mathcal R_3(v), \,\,\,\,\,\,\,\,\,\,Y_4(v) =L_0(v)+ \mathcal R_{1,4}, \text{ and }\,\,\,\,\,\, Y_5(v)=\mathcal R_2(v).
\ee
Due to the fast decay of $S_k$ and $\partial_tS_k$, the terms  $\mathcal R_2(v)$, $\mathcal R_3(v)$, and $\mathcal R_{1,4}$ are small (as stated in Proposition~\ref{apriori}. The term  $L_0(v)$ is not small but linear in $v$, and it is handled separately, see the proof of Theorem~\ref{thm:L1}. On the other hand one can not directly estimate the terms $Y_1$ and $Y_2$ without performing another differentiation by parts.   To do that we need to check the resonant terms in $Y_1$ and $Y_2$:
\begin{equation}\label{zeroset}
(k_1+k_2)(k_3+k_1)=0.
\end{equation}
The set for which \eqref{zeroset} holds is the disjoint union of the following 3 sets (taking into account that $k_2+k_3\neq 0$)
\begin{align*}
A_{1}&=\{k_1+k_2=0\}\cap\{k_3+k_1=0\}\Leftrightarrow \{k_1=-k,\ k_2=k,\ k_3=k\},\\
A_{2}&=\{k_1+k_2=0\} \cap \{ k_3+k_1\ne 0\} \Leftrightarrow \{k_1=j,\ k_2=-j,\ k_3=k,\ |j| \neq |k|\},\\
A_{3}&=\{k_3+k_1=0\}\cap\{k_1+k_2\ne 0\} \} \Leftrightarrow \{k_1=j,\ k_2=k,\ k_3=-j,\ |j| \neq |k|\}.
\end{align*}
We write
$$Y_2(v)_k=Y_{2r}(v)_k+Y_{2nr}(v)_{k}$$
where the subscript $r$ and $nr$ stands for the resonant and non-resonant terms respectively. We have
\begin{align*}
Y_{2r}(v)_k&=-\frac{i}{3}\sum_{\lambda=1}^{3}\sum_{A_{\lambda}}\frac{S_{k_1}v_{k_2}v_{k_3}}{k_1}=
-\frac{i}{3}\frac{S_{-k}v_{k}v_{k}}{-k} - \frac{2i}{3}v_k\sum_{\stackrel{j\in \Bbb Z_{0}}{ |j|\neq |k|}}\frac{S_{j}v_{-j}}{j}\\
& =\mathcal R_{1,2}(v)_k+ \mathcal R_{1,3}(v)_k,
\end{align*}
and
$$Y_{2nr}(v)_k=-\frac{i}{3}\sum_{k_1+k_2+k_3=k}^{nr}\frac{e^{-3it(k_1+k_2)(k_2+k_3)(k_3+k_1)}}{k_1}S_{k_1}v_{k_2}v_{k_3}.$$
Since the exponent in $Y_{2nr}(v)$ is not zero we can differentiate by parts one more time and obtain that
$$Y_{2nr}(v)_k =\partial_{t}M_{3}(v)_k+ M_{4}(v)_{k}$$
where
\be\label{m3}
M_{3}(v)_k= \frac{1}{9}\sum_{k_1+k_2+k_3=k}^{nr}\frac{e^{-3it(k_1+k_2)(k_2+k_3)(k_3+k_1)}}{k_1(k_1+k_2)(k_2+k_3)(k_3+k_1)}S_{k_1}v_{k_2}v_{k_3},
\ee
and
\be \nn\begin{split} M_{4}(v)_{k}= -\frac{1}{9}
\sum_{k_1+k_2+k_3=k}^{nr}\frac{e^{-3it(k_1+k_2)(k_2+k_3)(k_3+k_1)}}{k_1(k_1+k_2)(k_2+k_3)(k_3+k_1)}\times \\
\times \left( v_{k_2}v_{k_3}\partial_{t}S_{k_1}+S_{k_1}v_{k_3}\partial_{t}v_{k_2}+S_{k_1}v_{k_2}\partial_{t}v_{k_3}\right)
\end{split}
\ee
A calculation as before, by expressing time derivatives using \eqref{v_eq}, reveals that
\begin{align*}
M_4(v)_{k}&=-\frac{1}{9}\sum_{k_1+k_2+k_3=k}^{*}\frac{e^{-3it(k_1+k_2)(k_2+k_3)(k_3+k_1)}  v_{k_2}v_{k_3}\partial_{t}S_{k_1}}{k_1(k_1+k_2)(k_2+k_3)(k_3+k_1)},\\
&+\frac{i}{9}\sum_{k_1+k_2+k_3+k_4=k}^{\star}
\frac{e^{it\tilde\psi(k_1,k_2,k_3,k_4)}(k_3+k_4)S_{k_1}v_{k_2}(v_{k_3}+2S_{k_3})v_{k_4}}{k_1(k_1+k_2)(k_1+k_3+k_4)(k_2+k_3+k_4)}\\
&=\mathcal R_4(v)_k+\mathcal R_5(v)_k.
\end{align*}
The phase function $\tilde\psi$ is irrelevant for our calculations since it is going to be estimated out by taking absolute values inside the sums. For completeness we note that it can be expressed as
$$(k_1+k_2+k_3+k_4)^3-k_1^3-k_2^3-k_3^3-k_4^3.$$
Hence
\be\label{y2}
Y_{2}(v)=\mathcal R_{1,2}(v)_k+ \mathcal R_{1,3}(v)_k+\partial_{t}M_3(v) + \mathcal R_4(v)_{k}+\mathcal R_5(v)_{k}.
\ee
Similarly,
\be\label{y1}
Y_{1}(v)= \mathcal R_{1,1} +\partial_{t}N_3(v) +\mathcal R_6(v),
\ee
where
\be \label{n3}
N_{3}(v)_k:=\frac{1}{9}\sum_{k_1+k_2+k_3=k}^{*}\frac{e^{-3it(k_1+k_2)(k_2+k_3)(k_3+k_1)}v_{k_1}v_{k_2}v_{k_3}}{k_1(k_1+k_2)(k_2+k_3)(k_3+k_1)}.
\ee
Using \eqref{m3} and \eqref{n3}, note that $B_2(v)= -M_3(v) -N_3(v)$.
Therefore, substituting \eqref{y3y4y5}, \eqref{y2}, and \eqref{y1} in \eqref{y1y5},  we obtain \eqref{v_eq_dbp}.
\end{proof}
\begin{proof}[Proof of Proposition~\ref{dbpL1}]
We write the right hand side of \eqref{L1v} by distinguishing the resonant and nonresonant terms.
The resonant set corresponding to the terms with $k_2+k_3\neq 0$ is the same as above, and thus we get the following 3 terms:
$$
\frac{2i}{3}\Big(\frac{S_{-k}}{-k}S_kv_k+v_k\sum_{j,\,|j|\neq|k|}\frac{S_jS_{-j}}{j}+S_k\sum_{j,\,|j|\neq|k|}\frac{S_jv_{-j}}{j}\Big).
$$
Note that, by symmetry, the second term is zero. Combining the other terms we obtain the first summand in the definition of $E(v)$.
The terms with $k_2+k_3=0$ gives the second summand.

For the nonresonant terms we differentiate by parts as above obtaining
$$
 \frac{i}{3}\partial_t\Big(\sum_{k_1+k_2+k_3=k}^{*}\frac{e^{-3it(k_1+k_2)(k_2+k_3)(k_3+k_1)}}{k_1(k_1+k_2)(k_2+k_3)(k_3+k_1)} S_{k_1} S_{k_2}v_{k_3}\Big)
$$
$$-\frac{i}{3}\sum_{k_1+k_2+k_3=k}^{*}\frac{e^{-3it(k_1+k_2)(k_2+k_3)(k_3+k_1)}}{k_1(k_1+k_2)(k_2+k_3)(k_3+k_1)} \partial_t(S_{k_1} S_{k_2}v_{k_3}).
$$
The first line gives $D(v)$, and the second line gives the remaining terms in the definition of $E(v)$
after using the formula for $\partial_t v_{k_3}$ and renaming the variables.
\end{proof}

\section{Proof of Proposition~\ref{apriori}} \label{sec_apriori}
We start with the term $K$:
$$
\|K(v)\|_2\lesssim \|S_k/k\|_{\ell^1} \|v_k/k\|_2\lesssim \|v\|_{H^{-s}}.
$$
Similarly the $L^2$ norm of the first summand in the definition of $B(v)$ is
$$
\lesssim \|v_k/k\|_{\ell^1}\|v_k/k\|_2\lesssim \|v_k/|k|^s\|_2 \|1/|k|^{1-s}\|_{\ell^2}\|v_k/k\|_2 \lesssim \|v\|_{H^{-s}}^2,
$$
since $0<s<1/2$.
Note that the second summand is
\begin{align*}
&\lesssim    \sum_{k_1+k_2+k_3=k}^{*}\frac{|v_{k_1}|+|S_{k_1}|  }{|k_1||k_1+k_2||k_2+k_3||k_3+k_1|} |v_{k_2}||v_{k_3}| \\
&\lesssim    \sum_{k_1+k_2+k_3=k}^{*}\frac{ |v_{k_1}||v_{k_2}||v_{k_3}| }{|k_1||k_2|} +\frac{|k_1||S_{k_1}|v_{k_2}||v_{k_3}|}{|k_2||k_3|}.
\end{align*}
In the second line we use the following inequalities which are valid for the nonzero integral values of the factors:
$$|k_1+k_2||k_2+k_3||k_3+k_1|\gtrsim |k_2|,\,\,\,\,\, |k_1||k_1+n|\gtrsim |n|.$$
Therefore by Young's inequality, the $L^2$ norm of the second summand is
$$\lesssim \|v_k/k\|_{\ell^1}^2 \|v_k\|_{2}+\|kS_k\|_{\ell^1}\|v_k/k\|_{\ell^1}\|v_k/k\|_{2}\lesssim \|v\|_{H^{-s}}^2.$$
The $L^2$ bound for $L_0(v)$ follows from Young's inequality:
$$
\|L_0(v)\|_2\lesssim \|S_{k}/k\|_{\ell^1} \|S\|_{\ell^1} \|v\|_2\lesssim \|v\|_2.
$$
For the $H^{-s}$ bound, using the inequality $|k_3|^s\lesssim |k_1|^s|k_2|^s|k_1+k_2+k_3|^s=|k_1k_2k|^s$ (for nonzero integral values of the factors), we obtain
$$
\frac{|L_0(v)_k|}{|k^s|}\lesssim \sum_{k_1+k_2+k_3=k}\frac{|S_{k_1}|}{|k_1|^{1-s}} |k_2|^s |S_{k_2}| \frac{|v_{k_3}|}{|k_3|^s}.
$$
Therefore,
$$
\|L_0(v)\|_{H^{-s}}\lesssim \|S_k/|k|^{1-s}\|_{\ell^1}\||k|^s S_k \|_{\ell^1} \|v_k/|k|^s\|_2\lesssim \|v\|_{H^{-s}}.
$$

We now estimate $\mathcal R(v)$. Denote the terms in the $j$th line of the definition of $\mathcal R (v)$ by $\mathcal R_j (v)$, $j=1,2,...,6$.
The estimate for $\mathcal R_1 (v)$ follows as in the proof of Proposition~\ref{apriori2}.
The estimate for $\mathcal R_2 (v)$ is the same as the estimate for $K(v)$ with $S$ replaced with $\partial_t S$.
For $\mathcal R_3 (v)$, note that by Young's inequality
$$\|\mathcal R_3 (v) \|_2\lesssim \|v_k/k\|_{\ell^1} \|S\|_1 \|v\|_2\lesssim \|v\|_{H^{-s}}.$$
For $\mathcal R_4 (v)$, using $|k_1||k_1+k_2|\gtrsim |k_2|$, we have
$$\|\mathcal R_4 (v) \|_2\lesssim \|\partial_t S\|_{\ell^1} \|v_k/k\|_{\ell^1}  \|v\|_2\lesssim \|v\|_{H^{-s}}.$$
For $\mathcal R_5 (v)$, using $|k_1||k_1+k_2|\gtrsim |k_2|$, $|k_3+k_4|\lesssim |k_1||k_1+k_3+k_4|$ and $|k_1||k_2+k_3+k_4|\gtrsim |k_1+k_2+k_3+k_4|=|k|$, we have
$$|\mathcal R_5(v)_k|\lesssim \frac1{|k|}\sum_{k_1+k_2+k_3+k_4=k}^{\star} k_1^2 |S_{k_1}|\, \frac{|v_{k_2}|}{|k_2|}\, |v_{k_3}+2S_{k_3}|\, |v_{k_4}|.$$
Therefore, by Young's inequality
\begin{align*}
\|\mathcal R_5(v) \|_2& \lesssim \|1/k\|_{\ell^2} \Big\|\sum_{k_1+k_2+k_3+k_4=k}^{\star} k_1^2 |S_{k_1}|\, \frac{|v_{k_2}|}{|k_2|}\, |v_{k_3}+2S_{k_3}|\, |v_{k_4}|\Big\|_{\ell^\infty}\\
& \lesssim \| k^2 S_k\|_{\ell^1} \|v_k/k\|_{\ell^1}  \|v+2S\|_2 \|v\|_2 \lesssim \|v\|_{H^{-s}}.
\end{align*}
Finally, we consider $\mathcal R_6(v)$. Using $|k_1+2k_3+2k_4|\leq |k_1|+2|k_1+k_3+k_4|$, and then Cauchy Schwarz, we have
\begin{align*}
|\mathcal R_6(v)_k|^2 &\lesssim \Big[\sum_{k_1+k_2+k_3+k_4=k}^{\star}\frac{ |v_{k_1}v_{k_2}(v_{k_3}+2S_{k_3})v_{k_4}|}{|k_1+k_2||k_2+k_3+k_4|}
\Big(\frac1{|k_1|}+\frac1{|k_1+k_3+k_4|}\Big) \Big]^2\\
&\lesssim \sum_{k_1+k_2+k_3+k_4=k}^{\star}\frac{ v_{k_1}^2v_{k_2}^2(v_{k_3}+2S_{k_3})^2 }{|k_1|^{2s}}\times \\
& \times \sum_{k_1+k_2+k_3+k_4=k}^{\star}\frac {|k_1|^{2s}v_{k_4}^2}{|k_1+k_2|^2|k_2+k_3+k_4|^2}
\Big(\frac1{|k_1|^2}+\frac1{|k_1+k_3+k_4|^2}\Big)
\end{align*}
Note that the first factor above is $\lesssim \|v\|_{H^{-s}}^2 \|v+2S\|_2^2 \|v\|_2^2$. Therefore it suffices to prove that the sum of the second factor in
$k$ is $O(1)$. We write this sum as
\begin{align*}
& \sum_{k_1,k_2, k_3, k_4 }^{\star}\frac { v_{k_4}^2}{|k_1|^{2-2s}|k_1+k_2|^2|k_2+k_3+k_4|^2} \\
&  +\sum_{k_1,k_2,k_3,k_4}^{\star}\frac{|k_1|^{2s}v_{k_4}^2}{|k_1+k_2|^2|k_1+k_3+k_4|^2|k_2+k_3+k_4|^2}.
\end{align*}
To estimate the first sum, first take the sum in $k_3$, then in $k_2$, $k_1$, and $k_4$ in the given order.
The estimate for the second sum follows also by summing in the order given above and using the inequality
$$
\sum_{k_2,k_3}^{\star}\frac{1}{|k_1+k_2|^2|k_1+k_3+k_4|^2|k_2+k_3+k_4|^2}\lesssim \sum_{k_2}^*\frac{1}{|k_1+k_2|^2|k_1-k_2|^2}\lesssim \frac1{|k_1|^2}.
$$
Here we used the inequality (for $a\geq b > 1$)
$$
\sum_m^*\frac{1}{|n_1-m|^a|n_2-m|^b}\lesssim  \frac{1}{|n_1-n_2|^b}.
$$
\begin{rmk}\label{apriori_rho}
We note that in the proof of Proposition~\ref{apriori} and Proposition~\ref{apriori2}, the strongest conditions we need for $S_k=e^{-i\psi(k)t}\widehat{\Phi}_k$ are
$$\partial_t S_k\in\ell^1,\,\,\text{ and } k^2S_k\in \ell^1.$$
It is easy to see, using the equation, that $\partial_t \rho\in H^1$ for an $H^4$ solution $\rho$ of KdV.
Therefore, the sequence $e^{-i\psi(k)t}\widehat{\rho}_k$ satisfies the conditions above, and thus the assertions of  Proposition~\ref{apriori} and Proposition~\ref{apriori2} remains valid when $\Phi$ is replaced with $\rho(\cdot,t)\in H^4.$

\end{rmk}

{\bf Acknowledgement.} The authors would like to thank Michael Weinstein for a helpful discussion.

\end{document}